\input amstex
\input amssym

\magnification 1200
\loadmsbm
\parindent 0 cm

\define\nl{\bigskip\item{}}
\define\snl{\smallskip\item{}}
\define\inspr #1{\parindent=20pt\bigskip\bf\item{#1}}
\define\iinspr #1{\parindent=27pt\bigskip\bf\item{#1}}
\define\einspr{\parindent=0cm\bigskip}

\define\co{\Delta}

\define\ot{\otimes}

\define\<{\leg}
\define\>{\geq}

\define\om{\omega}

\define\s{\sum}
\define\si{\sigma}

\define\vp{\varphi}

\centerline{\bf The Larson-Sweedler theorem for}
\centerline{\bf multiplier Hopf algebras}
\bigskip
\bigskip
\centerline{\it Alfons Van Daele\rm($^{*}$) and \it Shuanhong Wang\rm($^{**}$)}
\bigskip
\bigskip
\bigskip
{\bf Abstract}
\bigskip
Any finite-dimensional Hopf algebra has a left and a right integral. Conversely, Larsen and Sweedler showed that, if a finite-dimensional algebra with identity and a comultiplication with counit has a faithful left integral, it has to be a Hopf algebra.
\smallskip
In this paper, we generalize this result to possibly infinite-dimensional algebras, with or without identity. We have to leave the setting of Hopf algebras and work with multiplier Hopf algebras. Moreover, whereas in the finite-dimensional case, there is a complete symmetry between the bialgebra and its dual, this is no longer the case in infinite dimensions. Therefore we consider a direct version (with integrals) and a dual version (with cointegrals) of the Larson-Sweedler theorem.
\smallskip
We also add some results about the antipode. Furthermore, in the process of this paper, we obtain a new approach to multiplier Hopf algebras with integrals.
\bigskip
\bigskip
August 2004 ({\it Version 1.0})
\bigskip
\vskip 6 cm
 \hrule
\medskip

($^{*}$) K.U.Leuven, Department of Mathematics, Celestijnenlaan
200B, BE-3001 Heverlee (Belgium). E-mail: Alfons.VanDaele\@wis.kuleuven.ac.be
\smallskip
($^{**})$ Southeast University,  Department of Mathematics, Nanjing 210096 (China). E-mail:
 shuanhwang2002\@yahoo.com, shuanhwang\@seu.edu.cn

\newpage

\bf 0. Introduction \rm
\snl 
Let $A$ be an algebra over $\Bbb C$ with identity $1$. Let
$\co: A\to A\ot A$ be a comultiplication on $A$. So $\co $ is a
unital homomorphism and satisfies coassociativity $(\co \ot \iota )\co =(\iota \ot \co )\co $
(where we use $\iota $ to denote the identity map). Let us also
assume the existence of a counit. It is a homomorphism $\varepsilon : A\to \Bbb C$
such that $(\varepsilon \ot \iota )\co (a)=a$ and $(\iota \ot \varepsilon )\co (a)=a$
for all $a\in A$.
\snl
Such a pair $(A, \co )$ is a Hopf algebra if there also exists an
antipode. This is a anti-homomorphism $S: A\to A$ such that
$m(S\ot \iota )\co (a)=\varepsilon(a)1$ and $m(\iota \ot S)\co (a)=\varepsilon(a)1$ for all
$a\in A$ (where we use $m$ to denote the multiplication map, defined
as a linear map from $A\ot A$ to $A$ by $m(a\ot b)=ab$).
\snl
It is well known that  any finite-dimensional Hopf algebra has a
left and a right integral (see [A], [S] and [VD4] for an alternative approach). Recall that a left integral is a non-zero linear
map $\varphi : A\to \Bbb C$ such that $(\iota \ot \varphi )\co(a)=\varphi(a)1$
for all $a\in A$ while a right integral is a non-zero linear map
$\psi : A\to \Bbb C$ such that $(\psi \ot \iota )\co (a)=\psi (a)1$
for all $a\in A$. These integrals are faithful in the sense that
the bilinear forms $(a, b)\mapsto \varphi (ab)$ and $(a, b)\mapsto \psi (ab)$
are non-degenerate (see e.g.\ Proposition 3.4 in [VD5] for a general proof of this fact).
\snl
Conversely, Larson and Sweedler showed in [L-S] that a pair $(A,
\co )$ of a finite-dimensional algebra with $1$ and a
comultiplication with a counit $\varepsilon$ is actually a Hopf algebra
if there exists a faithful left integral. Their result is stated
in its {\it dual form}, as we will explain in Section 4 of
this paper.
\snl
In this paper we treat two generalizations of this result. In the
first place, we no longer restrict to the finite-dimensional
case. The price we have to pay is that we need to assume the
existence of both a left and a right integral. Also because we no
longer assume the underlying algebra to be finite-dimensional,
there is a difference between the direct and the dual version.
The first one is about integrals (Section 2) while the second
one assumes the existence of cointegrals (Section 4). It also
seems to be more natural (and also more general) to leave
the framework of Hopf algebras (with unital algebras) and
pass to the more general theory of multiplier Hopf algebras (where the
algebras may not have an identity). We will also briefly mention the $^*$-algebra case but this is of minor importance for the treatment here.
\snl
The results, as well as the methods used to prove these results,
are greatly inspired by the recent developments in the theory
of locally compact quantum groups (see e.g.\ [K-V1] and [K-V2]). This theory can be considered
as the operator algebra approach to quantum groups. Before this theory was developed, there was
an obvious need for an operator algebra version of the notion of a
Hopf algebra. In operator algebras, people are
familiar with algebras without identity and the notion of a multiplier
algebra is commonly used. This inspired the first author
to generalize the notion of a Hopf algebra to algebras possibly
without identity and this lead to the theory of multiplier Hopf
algebras (see [VD2]).
\snl
But there is more. While the antipode in a purely algebraic
setting is not so hard to handle, it always has been a
troublesome object in the operator algebra approach. To begin
with, in interesting cases, the antipode becomes an unbounded
map, not everywhere defined. Moreover, basic formulas like
$m(S\ot \iota )\co (a)=\varepsilon(a)1$ become very hard to generalize
(for various reasons). This explains why, in the operator
algebra approach to quantum groups, usually the antipode is not
part of the axioms but rather an object that is constructed.
\snl
The most striking example of this we find in [K-V1, K-V2] where a theory
of locally compact quantum groups is developed. Roughly
speaking, a locally compact quantum group is defined as an
operator algebra with a comultiplication satisfying certain conditions.
One of these conditions is the existence of integrals. Then from
these axioms, the antipode is constructed.
\snl
An other example is found in [V-VD] where a theory of Hopf
C$^*$-algebras is developed. In this case, however, no integrals
are involved. Roughly speaking, this is the operator algebra
version of the theory of multiplier Hopf algebras, as developed
in [VD2]. Also in this theory, the antipode is not part of the
axioms but it is constructed.
\snl
The case of a locally compact quantum group is more like the
situation of the Larson-Sweedler theorem for finite-dimensional
Hopf algebras, but also the two other cases, where also the
antipode is constructed, have been sources of inspiration for this
paper. So, we have used some of the ideas and techniques from
these works on the operator algebra approach to quantum groups
in this paper. This again is a nice illustration of the advantage 
of having different approaches to quantum groups and the importance of the mutual influence.
\snl
Let us mention here the fact that the Larson-Sweedler result was first obtained by Kac and Palyutkin in the special case of a finite-dimensional operator algebra ([K-P]). See a remark in the introduction of [L-S].
\snl
We work in a purely algebraic setting however and no
knowledge about the operator algebra approach to quantum
groups is needed. So we will not use results from e.g [V-VD],
nor from [K-VD1, K-VD2]. We will just adopt some of the techniques to
this purely algebraic context. So, we mainly aim at algebraists
but it may be interesting for the reader to know
that we have some influence coming from analysis.
\nl
Let us now give a short survey of {\it the content of this paper}
and be a bit more precise about the main results and their
relation with the original work by Larson and Sweedler.
\snl
At the end of this introduction, we will give some basic
references and fix some notations used throughout this paper.
{\it Section 1} is devoted to some preliminaries. We recall the
definition of a comultiplication  $\co $ on an algebra $A$ possibly
without identity. We also recall the notion of a (regular)
multiplier Hopf algebra. New is the definition of the legs of $\co $
and the concept of a regular and of a full comultiplication.
\snl
In {\it Section 2} we prove the main result of the paper. We take an
algebra with a comultiplication. We assume the existence of integrals
and we show that we must have a multiplier Hopf algebra. We also
give another proof of (the dual form of) the original theorem of Larson and Sweedler
for a finite-dimensional algebra with a unit, a comultiplication
and a counit.
\snl
In {\it Section 3} we discuss properties of the counit and the antipode.
Recall that in the definition of a multiplier Hopf algebra, the
counit and the antipode are not present, they are constructed.
In this section we will give a new way to construct these
objects, using the integrals. As a byproduct, we obtain a
different approach to multiplier Hopf algebras with integrals
(sometimes called algebraic quantum groups, cf.\ [VD5]).
\snl
In {\it Section 4}, we treat the case of a cointegral. We again take an
algebra with a comultiplication. We first define the notion of a
left cointegral. We do not need a counit to define this concept.
In fact, our definition of a left cointegral makes it possible
to construct the counit. Then we assume  the existence of a left
and of a right cointegral and again we show that we have a
multiplier Hopf algebra (of discrete type, cf.\ [VD5, VD-Z1]). Also here we discuss
the relation of the antipode with the cointegrals. Specializing,
to the finite-dimensional case, we again find a proof of (the original form of) the
Larson-Sweedler theorem as it is found in [L-S].
\nl
We finish this introduction by collecting some of the {\it basic
notions} for this paper and give some standard references.
\snl
We work with algebras over $\Bbb C$ with or without identity.
If the algebra has an identity, we denote it by $1$. If there is
no identity, we want the product to be non-degenerate, i.e.\ if $x\in A$
then $x=0$ if either $ax=0$ for all $a\in A$ or $xa=0$ for all $a\in
A$. The product in an algebra with identity is automatically
non-degenerate. For any algebra $A$ we have the multiplier
algebra $M(A)$. It can be characterized as the largest algebra with
identity such that $A$ sits in $M(A)$ as an essential two-sided
ideal. By {\it essential} we mean that if $x\in M(A)$ then $x=0$ if
either $xa=0$ for all $a\in A$ or $ax=0$ for all
$a\in A$. If $A$ has an identity, then $M(A)=A$.
\snl
If $A$ and $B$ are algebras as above and $\gamma : A\to M(B)$ is an
algebra map, then $\gamma $ is called non-degenerate if $\gamma (A)B=B\gamma
(A)=B$. If this is the case, then $\gamma $ has a unique extension to
a unital homomorphism from $M(A)$ to $M(B)$. This extension is
still denoted by $\gamma $. For details, see e.g.\ the appendix in [VD2].
\snl
We will use the tensor product $A\ot A$ for an algebra $A$. It is
again an algebra with a non-degenerate product. So we can also
define the multiplier algebra $M(A\ot A)$. A comultiplication on
$A$ is a homomorphism from $A$ to $M(A\ot A)$.
All the comultiplications we encounter in this paper will turn
out to be non-degenerate. So they can be extended to unital
maps from $M(A)$ to $M(A\ot A)$. Then coassociativity can be stated
in the form $(\co \ot \iota )\co =(\iota \ot \co )\co $ where $\iota $ is the
identity map and $\co \ot \iota $ and $\iota \ot \co$ are the extensions
of the corresponding homomorphisms from $A\ot A$ to $M(A\ot A\ot
A)$. There is also another way to look at coassociativity
(see Definition 1.1 in Section 1 of this paper).
\snl
Occasionally, we will use the Sweedler notation $\co (a)=\sum a_{(1)}\ot
a_{(2)}$ and $(\co \ot \iota )\co (a)=\sum a_{(1)}\ot a_{(2)}\ot a_{(3)}$.
The use of the Sweedler notation in the case of multiplier Hopf algebras 
has been justified e.g.\ in [D-VD] but it has to be done with
some care.
\nl
Let us now collect some {\it standard references}. For the general theory of Hopf algebras, we refer to [A] and [S]. The original reference for multiplier Hopf algebras is [VD2] but a survey can be found in [VD-Z2]. Integrals on multiplier Hopf algebras have been studied in [VD5] and they are also considered in the survey paper [VD-Z2]. Multiplier Hopf algebras with cointegrals (i.e.\ multiplier Hopf algebras of discrete type) where already introduced in [VD5] but have been studied in [VD-Z1]. 
\nl
\nl
{\bf Acknowlegdment} The first author would like to thank the organisers of the meeting {\it Groupes quantiques localement compacts et sym\'etries quantiques} (June 2004 at CIRM, Marseille) for the invitation and the opportunity to give a talk related to this work. The second author was supported by the Research Council of the K.U.Leuven.
\nl\nl

\bf 1. Preliminaries \rm
\nl
Let $A$ be an algebra over $\Bbb C$, with or without unit,
but with a  non-degenerate product. Consider the tensor product $A\ot A$
of $A$ with itself. The product in $A\ot A$ will again be
non-degenerate. We can consider the multiplier algebras $M(A)$
and $M(A\ot A)$ of $A$ and $A\ot A$ respectively. There are natural
imbeddings of $A\ot A$ in $M(A)\ot M(A)$ and of $M(A)\ot M(A)$ in
$M(A\ot A)$. In general, these two imbeddings are strict.
\snl
This is the setting we need to define what we will mean by a
(regular) comultiplication in this paper (see e.g.\ [VD2]):

\inspr{1.1}  Definition \rm 
A homomorphism $\co : A\to
 M(A\ot A)$ is called a {\it comultiplication} if
\snl
i)\quad  $\co(a)(1\ot b)\in A\ot A$ and $(a\ot 1)\co(b)\in A\ot A$ for
 all $a, b\in A$,
\snl
ii)\quad $(a\ot 1\ot 1)(\co \ot \iota )(\co(b)(1\ot c))
   =(\iota \ot \co )((a\ot 1)\co (b))(1\ot 1\ot c)$ for all $a, b, c\in A$ ({\it coassociativity}).
\einspr

Observe that i) is needed to give a meaning to ii). As we
mentioned already in the introduction, we use $1$ for the
identity (here in $M(A)$) and $\iota $ for the identity map
(here from $A$ to itself). If $A$ has an identity, then
i) is automatic and ii) is nothing else but coassociativity
$(\co \ot \iota )\co =(\iota \ot \co )\co $. 
\snl
As we have mentioned in the introduction, occasionally, we will say something about the $^*$-algebra case. When we have a comultiplication $\co$ on a $^*$-algebra, we always require $\co$ to be a $^*$-homomorphism.
\snl
In this paper, we will only work with regular comultiplications:

\inspr{1.2}  Definition \rm 
A comultiplication $\co $ on $A$ is called {\it regular} if also
\snl
iii) \quad $\co (a)(b\ot 1)\in A\ot A$
 and $(1\ot a)\co (b)\in A\ot A$ for all $a, b\in A$.
\einspr

For an algebra with identity, regularity is automatic. Also in the case of a $^*$-algebra, a comultiplication is automatically regular in the above sense.
\snl
We have the following result.

\inspr{1.3}  Proposition \rm 
If $\co $ is a regular comultiplication,
then the opposite comultiplication $\co' $, obtained by composing
$\co $ with the flip map, is again a comultiplication.

\snl\bf Proof: \rm 
First recall the definition of the opposite
comultiplication. The flip map $\si : A\ot A\to A\ot A$ is
defined as usual by $\si (a\ot b)=b\ot a$. It can be extended
uniquely to a homomorphism, still denoted by $\si $,
from $M(A\ot A)$ to itself. Then $\co'$ is defined by
$\co'(a)=\si (\co (a))$ for all $a\in A$.
\snl
Now, it is obvious that condition iii) in Definition 1.2 gives
condition i) in Definition 1.1 for $\co'$. To prove that
$\co'$ is also coassociative, one simply has to start with
the formula in ii), multiply with $1\ot 1\ot c'$ on the left and
with $a'\ot 1\ot 1$ on the right, bring where possible $1\ot c'$
and $a'\ot 1$ inside the brackets, take $a\ot 1$ and $1\ot c$
outside and finally use the non-degeneracy of the product and
cancel $a\ot 1\ot 1$ on the left and $1\ot 1\ot c$ on the right. This
will yield condition ii) for $\co'$ provided we flip
the first and the third factor. 
\einspr

Let us now recall the notion of a (regular) multiplier Hopf
algebra (see [VD2]).

\inspr{1.4}  Definition \rm 
Let $A$ be an algebra over
$\Bbb C$ with a non-degenerate product and let $\co $ be
a comultiplication on $A$. Then $(A, \co )$ is called a
{\it multiplier Hopf algebra} if the linear maps $T_1, T_2$,
defined from $A\ot A$ to $A\ot A$ by
$$ T_1(a\ot b)=\co (a)(1\ot b) \qquad\text{and}\qquad T_2(a\ot b)=(a\ot 1)\co (b) $$
are bijective. A multiplier Hopf algebra is called {\it regular}
if $\co $ is a regular comultiplication and if also the maps
above, now for $\co'$, are bijections.
\einspr

Remark that any Hopf algebra is a multiplier Hopf algebra
and conversely, any multiplier Hopf algebra,
with an algebra with identity, is a Hopf algebra (see [VD2]).
Regularity is not automatic for Hopf algebras, but equivalent
with the antipode being bijective. This is automatic in the case of a multiplier Hopf $^*$-algebra (i.e.\ when the underlying algebra is a $^*$-algebra).
\snl
We will need the notion of {\it the legs} of $\co $:

\inspr{1.5}  Definition \rm 
Let $A$ be an algebra with a non-degenerate
product and $\co $ a comultiplication on $A$. For
any $a\in A$ we define {\it the left leg}  of $\co (a)$ as
the smallest subspace  $V$ of $A$ so that
$$\co (a)(1\ot A)\subseteq V\ot A.$$
Similarly, {\it the right leg}  of $\co (a)$ is
the smallest subspace  $W$ of $A$ so that
$$(A\ot 1)\co (a)\subseteq A\ot W.$$
Similarly, we define the left and the right legs of $\co $
as the smallest subspaces $V$ and $W$ of $A$ satisfying
$$\co (A)(1\ot A)\subseteq V\ot A \qquad\text{and}\qquad (A\ot 1)\co (A)\subseteq A\ot W.$$
\einspr

In the case of a multiplier Hopf algebra, we clearly have that
the left and the right leg of $\co $ are all of $A$ (because the maps
$T_1$ and $T_2$ are surjective). In fact using the counit $\varepsilon$,
we see that $a$ belongs to both the right and the left leg of $\co (a)$.
\snl
In general, we can show the following.

\inspr{1.6}  Proposition \rm 
The left leg of $\co(a)$ is spanned by
elements of the form $(\iota \ot \om)(\co (a)(1\ot b))$ where
$b\in A$ and $\om \in A'$ (the dual space of $A$).
Similarly, the right leg of $\co (a)$ is spanned by the elements
$(\om \ot \iota )((b\ot 1)\co (a))$ where $b\in A$ and $\om \in A'$.

\snl \bf Proof: \rm 
It is clear that $(\iota \ot \om)(\co (a)(1\ot b))$ belongs
to the left leg of $\co (a)$ for all $b\in A$ and $\om \in A'$.
Conversely, let $V$ be the space spanned by these elements. We
need to show that $\co (a)(1\ot b)\in V\ot A$ for all $b\in A$.
Take $b\in A$ and write $\co (a)(1\ot b)=\sum p_i\ot q_i$
with the $\{ q_i\}$ linearly independent. Choose
$\om \in A'$ such that $\om (q_i)=1$ for some $i$ and $\om (q_j)=0$
for $j\ne i$. Then $p_i=(\iota \ot \om )(\co (a)(1\ot b))$. By
assumption $p_i\in V$. This is true for all $i$ and so
$\co (a)(1\ot b)\subseteq V\ot A$.
\snl
Similarly for the right leg of $\co (a)$. 
\einspr

For a regular comultiplication, we have the following.

\inspr{1.7}  Proposition \rm 
If $\co $ is  a regular comultiplication
and $a\in A$, then the left leg of $\co (a)$ is also the
smallest subspace $V'$ of $A$ so that
$$(1\ot A)\co (a)\subseteq V'\ot A.$$

\snl\bf Proof: \rm 
Suppose that $V'$ is a subspace of $A$ so that $(1\ot
A)\co(a)\subseteq V'\ot A$. We will show that $\co (a)(1\ot b)\in V'\ot
A$ for all $b\in A$. Then it will follow that the left leg of $\co (a)$
is contained in $V'$. A similar argument will give that $V'$ is
contained in the left leg of $\co (a)$ and this will
prove the result.
\snl
So, let $b\in A$ and write $\co (a)(1\ot b)=\sum p_i\ot q_i$ with
the $\{q_i\}$ independent. We know that $\sum p_i\ot cq_i=
(1\ot c)\co (a)(1\ot b)\in V'\ot A$ for all $c\in A$.
Assume that $\om \in A'$ and that $\om (x)=0$ for all
$x\in V'$. Then $\sum \om (p_i)cq_i=0$ for all $c$.
By the non-degeneracy of the product, we get
$\sum \om (p_i)q_i=0$.  As the $\{q_i\}$ are chosen to be
linearly independent, it follows that $\om (p_i)=0$ for all $i$.
This is true for any $\om $ that vanishes on $V'$. Therefore
$p_i\in V'$ for all $i$. Hence $\co (a)(1\ot b)\in V'\ot A$. This
proves the claim. 
\einspr

Of course, we have a similar result for the right leg of $\co(a)$
and similar results as in Proposition 1.6, now with elements
$(\iota \ot \om )((1\ot b)\co (a))$ for the left leg and
$(\om \ot \iota )(\co (a)(b\ot 1))$ for the right leg. We also have
similar results for the left and the right leg of $\co $. 
\snl
If $A$ is a $^*$-algebra, then the comultiplication is automatically regular so that Proposition 1.7 applies. Then, if $a$ is a self-adjoint element, 
i.e.\ $a=a^*$, the left and right legs of $\co(a)$ are self-adjoint subspaces of $A$.

\nl
We finish
this preliminary section with the following definition.

\inspr{1.8}  Definition \rm 
Let $\co $ be a comultiplication on $A$. We
call it {\it full} if the left and the right legs of $\co $ are all of
$A$.
\einspr

As we have seen already, when $(A, \co )$ is a multiplier Hopf
algebra, then $\co $ is full. In fact this will already be the
case when there is a counit (i.e.\ a linear map $\varepsilon$ so that
$(\iota \ot \varepsilon )\co (a)=a$, in the sense that $(\iota \ot \varepsilon )((b\ot 1)\co (a))=ba$
for all $a, b\in A$ and so that $(\varepsilon \ot \iota )\co (a)=a$
in the sense that $(\varepsilon \ot \iota )(\co (a)(1\ot b))=ab$).
\snl
We also have a full comultiplication if the maps $T_1$ and $T_2$,
defined as before by
$$T_1(a\ot b)=\co (a)(1\ot b) \qquad\text{and}\qquad T_2(a\ot b)=(a\ot 1)\co (b)$$
are surjective. However, it should be noticed that, in general,
the surjectivity of these maps is not implied when $\co $ is full.
In the next section, we will see that this is the case when
there exist integrals.
\nl\nl

\bf 2. The Larson-Sweedler theorem \rm
\nl
Throughout this section, $A$ will be an algebra over $\Bbb C$,
with a non-degenerate product and $\co $ will be a
comultiplication on $A$. Most of the time, we will assume
that $\co $ is regular.
\snl
Let us recall the definition of an integral in this setting.

\inspr{2.1} Definition \rm 
A linear functional $\vp $ on $A$ is called
{\it left invariant} if
$$(\iota \ot \vp )((b\ot 1)\co (a))=\vp (a)b$$
for all $a, b\in A$. A non-zero left invariant functional is
called a {\it left integral}. Similarly, a linear functional $\psi $
is called  {\it right invariant} if
$$(\psi \ot \iota )(\co (a)(1\ot b))=\psi (a)b$$
for all $a, b\in A$. A non-zero right invariant functional is a {\it right integral}.
\einspr

In the case of a $^*$-algebra, it is common to assume positivity of these integrals. This means that $\vp(a^*a)\geq 0$ for all $a$ and similarly for $\psi$. In any case, one can always assume that $\vp(a^*)=\overline{\vp(a)}$ without loss of generality.
\snl
In the case of a regular comultiplication, for every $a\in A$ and
$\om \in A'$ (the dual space of $A$), we can define elements $(\iota \ot \om )\co (a)$ and
$(\om \ot \iota )\co (a)$ in $M(A)$ in the obvious way. Then
$\vp \in A'$ is left invariant if $(\iota \ot \vp )\co(a)=\vp (a)1$
for all $a\in A$ and $\psi \in A'$ is right invariant if
$(\psi \ot \iota )\co (a)=\psi (a)1$ for all $a\in A$.
\snl
We will also use the following terminology.

\inspr{2.2}  Definition \rm
A linear functional $\om $ on $A$ is called
{\it faithful} if the bilinear map $(a, b)\mapsto \om (ab)$ is
non-degenerate.
\einspr

So $\om $ is faithful if $b=0$ when $\om (ab)=0$ for all $a$ or
when $\om (ba)=0$ for all $a$.
\snl
In the case of a finite-dimensional algebra, we only need one
of these conditions. Indeed, suppose that $b\mapsto \om (\,\cdot\, b)$ is
injective from $A$ to $A'$. Then it is also surjective. If
now $\om (ab)=0$ for all $b\in A$, then $\rho(a)=0$ for all $\rho\in A'$
and so $a=0$.
\snl
In the case of a $^*$-algebra and a positive linear functional $\om$, faithfulness means that $\om(a^*a)=0$ implies $a=0$.
\snl
It was shown in [VD5, Proposition 3.4] that integrals on multiplier Hopf algebras are automatically faithfull. They are also unique [VD5, Theorem 3.7].
\snl
Now we are able to formulate the main result.

\inspr{2.3}  Theorem \rm 
Let $A$ be an algebra with a non-degenerate
product and assume that $\co $ is a full and regular
comultiplication  on $A$. If there exists a faithful left
integral and a faithful right integral, then
$(A, \co )$ is a regular multiplier Hopf algebra.
\einspr
Later in this section, we will look at the finite-dimensional case and discuss this result. In particular, we will relate it with the original theorem of Larson and Sweedler. The reader who is interested in the deeper relation with the operator algebra approach to quantum groups, is invited to compare the assumptions in this theorem with the axioms of a locally compact quantum group in [K-V1] and [K-V2]. 
\snl
Now, we start with the proof of the theorem. We split up the proof in a few lemmas as we will need
these smaller results later when we discuss the theorem in the finite-dimensional case.

\inspr{2.4} Lemma \rm 
Let $\co $ be any comultiplication on $A$. Assume
that there is a faithful right integral. Then the linear map
$T_1$, defined by $T_1(a\ot b)=\co (a)(1\ot b)$,  is injective.

\snl\bf Proof: \rm 
Suppose $\sum \co (a_i)(1\ot b_i)=0$. Take any $x\in A$ and multiply with
$\co (x)$ from the left and apply $\psi \ot \iota $ where
$\psi $ is a right integral. We get, using the invariance, that
$\sum \psi (xa_i)b_i=0$. If now $\psi $ is faithful, we get
$\sum a_i\ot b_i=0$. This proves the injectivity of $T_1$.
\einspr

In a similar way, when $\co $ is regular, we get the injectivity
of the map $a\ot b\mapsto (1\ot a)\co (b)$ when we have a faithful
right integral. When we have a faithful left integral, we get
$T_2$ injective with $T_2(a\ot b)=(a\ot 1)\co (b)$. And when $\co $
is regular, also $a\ot b\mapsto \co (a)(b\ot 1)$ will be injective.
\snl
In other words we get:

\inspr{2.5}  Proposition \rm 
If $\co $ is a regular comultiplication and
if there exist faithful left and right integrals, then the
four maps
$$\matrix
&a\ot b\mapsto \co (a)(1\ot b)& \quad \quad 
&a\ot b\mapsto \co (a)(b\ot 1)\\
&a\ot b\mapsto (a\ot 1)\co (b)& \quad \quad 
&a\ot b\mapsto (1\ot a)\co(b)
\endmatrix$$
are all injective.
\einspr

In the $^*$-algebra case, the injectivity of the maps on the left will already imply the injectivity of the maps on the right.
\snl
What about surjectivity? We get the following.

\inspr{2.6} Lemma \rm 
Let $\co $ be a regular comultiplication
and let $\vp $ be a left integral. Take $a, b\in A$ and
define $x=(\iota \ot \vp )(\co (a)(1\ot b))$. Then
$x\ot c$ belongs to the range of $T_1$ for all $c\in A$.

\snl\bf Proof: \rm 
We claim that
$x\ot c=T_1(y)$
where
$$y=(\iota \ot \iota \ot \vp )(\co _{13}(a)\co _{23}(b)(1\ot c\ot 1))$$
with $\co_{23}(b)=1\ot \co (b)$ and $\co_{13}(a)=(1\ot \si )(\co(a)\ot 1)$ where $\si $ is the flip map.
\snl
First let us verify that $y$ is well-defined in $A\ot A$.
Because $\co $ is assumed to be regular, we have $\co (b)(c\ot 1)\in A\ot
A$. For all $p, q\in A$ we get $\co_{13}(a)(1\ot p\ot q)\in A\ot A\ot
A$. So $\co _{13}(a)\co _{23}(b)(1\ot c\ot 1)\in A\ot A\ot A$ and if we
apply $\iota \ot \iota \ot \vp $, we get $y\in A\ot A$.
\snl
Now
$$(T_1\ot \iota )(\co_{13}(a)(1\ot p\ot q))=
((\co \ot \iota )\co (a))(1\ot p\ot q)
$$
and so
$$\align
(T_1\ot \iota )(\co_{13}(a)\co_{23}(b)(1\ot c\ot 1))
&=((\co \ot \iota )\co (a))((1\ot \co (b))(1\ot c\ot 1))\\
&=(\iota \ot \co )(\co (a)(1\ot b))(1\ot c\ot 1).
\endalign$$
If we apply $\iota \ot \iota \ot \vp $, it will follow from the left
invariance of $\vp $ that
$$T_1(y)=(\iota \ot \vp )(\co (a)(1\ot b))\ot c=x\ot c$$
and this proves the result. 
\einspr

Now we use this lemma to prove that $T_1$ is surjective under
certain conditions.

\inspr{2.7}  Lemma \rm
Let $\co $ be a regular comultiplication
such that the left leg of $\co $ is all of $A$. If there
is a faithful left integral, then $T_1$ is surjective.

\snl\bf Proof: \rm 
Because of the previous lemma, we have to show that
$A$ is spanned by the elements $(\iota \ot \vp )(\co (a)(1\ot b))$
where $a, b\in A$. Suppose that $\om \in A'$ and
that $\om ((\iota \ot \vp )(\co (a)(1\ot b)))=0$ for all $a, b\in A$.
Then $\vp ((\om \ot 1)(\co (a)(1\ot b))\, b')=0$
for all $a, b, b'\in A$. Because $\vp $ is faithful,
$(\om \ot 1)(\co (a)(1\ot b))=0$ for all $a, b\in A$.
Then $\om $ is $0$ on the left leg of $\co $. By assumption
$\om =0$. It follows that indeed $A$ is spanned by
elements of the form $(\iota \ot \vp )(\co (a)(1\ot b))$
with $a, b\in A$. This proves the lemma. 
\einspr

In a completely similar fashion, we will get that the
following is true.

\inspr{2.8}  Proposition \rm 
Let $\co $ be a regular and full comultiplication. Assume that there exist a faithful left integral
and a faithful right integral. Then the four maps
$$\matrix
&a\ot b\mapsto \co (a)(1\ot b) &\quad \quad 
&a\ot b\mapsto \co (a)(b\ot 1)\\
&a\ot b\mapsto (a\ot 1)\co (b) &\quad \quad 
&a\ot b\mapsto (1\ot a)\co(b)
\endmatrix$$
are surjective.
\einspr

Again, in the $^*$-algebra case, the surjectivity of the maps on the left will imply already the surjectivity of the maps on the right.
\snl

Now, we have also completed the proof of the theorem. Indeed, if
$\co $ is regular and full and if we have a faithful left integral
and a faithful right integral, all these four maps are bijective.
Then, $(A, \co )$ is a regular multiplier Hopf algebra (see [VD2]).
\snl
In the next section, we will see how the antipode and the counit
can be obtained in this case, without using the general theory.
\nl
Now, we finish this section by looking at the finite-dimensional
case.

\inspr{2.9}  Theorem \rm 
Let $A$ be a  finite-dimensional algebra
with $1$ and $\co $ a comultiplication on $A$. Assume that the left leg of $\co $ is all of $A$.
If there is a faithful left integral, then
$(A, \co )$ is a Hopf algebra.

\snl\bf Proof: \rm
By the analogue of Lemma 2.4, the map
$T_2$, defined by $T_2(a\ot b)=(a\ot 1)\co (b)$ is injective.
Because we are  working in finite dimensions, this map is also
surjective. By Lemma 2.7 the map $T_1$, defined by
$T_1(a\ot b)=\co (a)(1\ot b)$ is surjective. Again here this
implies that it is also injective. Hence
$(A, \co )$ is a Hopf algebra.  
\einspr

If we have a finite-dimensional algebra $A$ with $1$ and
a comultiplication $\co $ with a counit $\varepsilon$, then this
comultiplication is automatically full. So in this
case the extra assumption in the theorem above (about the left
leg of $\co $) is fulfilled. Then we get a (dual) version
of the original Larson-Sweedler theorem, as we find it in [L-S].
The faithfulness of the left integral in our result corresponds with the non-singularity assumption in their theorem. 
In Section 4 we will recover the Larson-Sweedler as it is
formulated in the original paper (and we will give some more comments).
\nl\nl

\bf 3. The counit and the antipode \rm
\nl
In [VD2], we have constructed the counit and the antipode for a
multiplier Hopf algebra. Let us recall how this was done (see
Section 3 and 4 in [VD2]).

\inspr{3.1}  Proposition \rm 
Let $(A, \co )$ be a multiplier Hopf algebra.
There exists a linear map $\varepsilon: A\to \Bbb C$ defined by
$\varepsilon(a)b=\sum p_iq_i$ when $a\ot b=\sum \co (p_i)(1\ot q_i)$. There
exists a linear map $S: A\to M(A)$ defined by
$S(a)b=\sum \varepsilon (p_i)q_i$ when $a\ot b=\sum \co (p_i)(1\ot q_i)$.
\einspr

The bijectivity of the map $T_1$, as defined already in Section 2
by $T_1(a\ot b)=\co (a)(1\ot b)$, is one of the assumptions of a
multiplier Hopf algebra and is used in the proof of the
proposition above. If the multiplier Hopf algebra is regular, it is shown that $S$ maps $A$ into $A$ and that it is bijective.
\snl
In the case where we have a regular multiplier Hopf algebra with
a left integral, we have the following result.

\inspr{3.2}  Proposition \rm 
If $(A, \co )$ is a regular multiplier Hopf
algebra and $\vp $ a left integral, then
$$\varepsilon ((\iota \ot \vp )(\co (a)(1\ot b)))=\vp (ab)$$
and
$$ S((\iota \ot \vp )(\co (a)(1\ot b)))=(\iota \ot \vp )((1\ot a)\co (b))$$
for all $a, b\in A$.
\einspr

The first formula follows trivially from the property of the
counit saying
$$(\varepsilon \ot \iota )(\co (a)(1\ot b))=ab.$$
Also the second formula can be proven from the standard
properties of the antipode.
\snl
However, here let us use the formula obtained in Lemma 2.6. We
showed that
$$x\ot c=T_1((\iota \ot \iota \ot \vp )(\co_{13}(a)\co_{23}(b)(1\ot c\ot 1)))
$$
where $x=(\iota \ot \vp )(\co (a)(1\ot b))$. Then it follows
already from the defining formulas in Proposition 3.1 that
$$
\varepsilon (x)c=(\iota \ot \vp )(\co (a)\co (b)(c\ot 1))
=\vp (ab)c,$$
so that $\varepsilon(x)=\vp (ab)$, and
$$\align
S(x)c
&=(\varepsilon \ot \iota \ot \vp )(\co_{13}(a)\co_{23}(b)(1\ot c\ot 1))\\
&=(\iota \ot \vp )((1\ot a)\co (b)(c\ot 1))
\endalign$$
so that $S(x)=(\iota \ot \vp )((1\ot a)\co (b))$.
\nl
In the remaining of this section, as we announced already in the introduction, we will use the formulas in
Proposition 3.2 to construct the counit and the antipode when we
have a pair of an algebra $A$ with a coproduct as in Theorem 2.3 and to prove their properties. We do not get any new results but it seems nice and instructive.
We do not use any of the results from [VD5]. We will mainly take advantage of the extra information we get from the integrals. So, we get an alternative approach to the theory of multiplier Hopf algebras with integrals (the so-called algebraic quantum groups). Whenever convenient, we feel free to use the classical arguments.
\nl
So, as in Theorem 2.3, let $A$ be an algebra over $\Bbb C$,
with a non-degenerate product, and let $\co $ be a regular
comultiplication which is full, that is such that the two legs
of $\co $ are all of $A$. Also assume that we have a faithful left
integral $\vp $ and a faithful right integral $\psi $. We know
that these integrals must be unique (up to a scalar), but we will
not use this information.
\snl
First, we treat the counit. To define the counit, we need the
following lemma.

\inspr{3.3}  Lemma \rm
If $\sum (\iota \ot \vp )(\co (a_i)(1\ot b_i))=0$, then
also $\sum \vp (a_ib_i)=0$.

\snl\bf Proof: \rm
We know that for all $a, b, c\in A$ we have
$$x\ot c = 
T_1((\iota \ot \iota \ot \vp )(\co_{13}(a)\co_{23}(b)(1\ot c\ot 1)))$$
where $x=(\iota \ot \vp )(\co (a)(1\ot b))$. From the existence of
a faithful right integral, we know that $T_1$ is injective.
Finally observe that
$$\align
m(\iota \ot \iota \ot \vp )(\co_{13}(a)\co_{23}(b)(1\ot c\ot 1))
&=(\iota \ot \vp )(\co (a)\co (b)(c\ot 1))\\
&=\vp (ab)c.
\endalign$$
If we apply all this to $a_i, b_i$ and take the sum, we get the
result. 
\einspr

Then, the following definition makes sense.

\inspr{3.4}  Definition \rm 
Define a linear map $\varepsilon : A\to \Bbb C$
by
$$
\varepsilon ((\iota \ot \vp )(\co (a)(1\ot b)))=\vp (ab)
$$
whenever $a, b\in A$.
\einspr

Recall that every element in $A$ is a sum of elements of the form
$(\iota \ot \vp )(\co (a)(1\ot b))$ with $a, b\in A$ (see the proof of
Lemma 2.7). So, with this definition, the counit is indeed defined on all of $A$.
\snl
Also notice that, from the proof of the lemma, we get with this definition that indeed $\varepsilon(x)c=\sum p_iq_i$ if $x\ot c=\sum\co(p_i)(1\ot q_i)$. This is in accordance with the way $\varepsilon$ is originally introduced (see Proposition 3.1). It also shows that the above definition does not depend on the choice of $\vp$.
\snl
It also follows immediately from the definition and
the faithfulness of $\vp $ that $(\varepsilon \ot \iota )\co (a)=a$ for all
$a\in A$. Recall that, as we saw in a remark after Definition 2.1, elements of the form $(\omega\ot\iota)\co(a)$ and $(\iota\ot\omega)\co(a)$ are well defined in $M(A)$ for all $a\in A$ and any linear functional  $\omega$ on $A$. In particular, this is the case with $\omega=\varepsilon$.
\snl
Let us now try to prove the other main properties of
the counit.

\inspr{3.5}  Proposition \rm
The map $\varepsilon : A\to \Bbb C$ is a
homomorphism and satisfies
$$(\varepsilon \ot \iota )\co (a)=a \qquad\text{and}\qquad
 (\iota \ot \varepsilon )\co (a)=a$$
for all $a$. Any other linear map satisfying one of these two
equations for all $a$ must be $\varepsilon $.

\snl\bf Proof: \rm 
There are different ways to obtain that $\varepsilon $ is a
homomorphism. A cute way is as follows. It is inspired by techniques used in [V-VD]. 
\snl
Take $a, b, c$ in $A$
and use the surjectivity of $T_1$ to write
$$\align
ab\ot c = (a\ot c)(b\ot 1)
&=(\sum \co (p_i)(1\ot q_i))(b\ot 1)\\
&=\sum \co (p_i)(b\ot q_i)\\
&=\sum \co (p_i)\co (r_{ij})(1\ot s_{ij})\\
&=\sum \co (p_ir_{ij})(1\ot s_{ij}).
\endalign$$
Then, by the remark following Definition 3.4, we get
$$\align
\varepsilon(ab)c &= \sum p_ir_{ij}s_{ij}
=\sum p_i\varepsilon (b)q_{i}\\
&=\varepsilon (b)\sum  p_iq_i=\varepsilon (b)\varepsilon (a)c.
\endalign$$
This shows that $\varepsilon $ is a homomorphism.
\snl
We have seen already that $(\varepsilon \ot \iota )\co (a)=a$. Because of the
symmetry of our data, we also have a linear map $\varepsilon': A\to \Bbb C$
satisfying  $(\iota \ot \varepsilon')\co (a)=a$ for all $a$. Now, suppose that $\om $ is any linear map from $A$ to $\Bbb C$ satisfying
$(\iota \ot \om )\co (a)=a$. Then we have for all $a, b, c\in A$ that
$$\align
(\varepsilon \ot \iota )((c\ot 1)\co (a)(1\ot b)) 
&=\varepsilon (c)(\varepsilon \ot \iota )(\co (a)(1\ot b))\\
&=\varepsilon (c)ab.
\endalign$$
We can cancel $b$ and obtain
$$(\varepsilon \ot \iota )((c\ot 1)\co (a))=\varepsilon (c)a$$
for all $a, c\in A$. If we apply $\om $ we get
$$\varepsilon (ca)=\varepsilon (c)\om (a)$$
for all $a, c\in A$. As $\varepsilon (ca)=\varepsilon (c)\varepsilon (a)$ and
$\varepsilon$ is non-zero, we get $\varepsilon (a)=\om (a)$. This
not only proves the uniqueness statement but also
that $(\iota \ot \varepsilon )\co (a)=a$ because $\varepsilon =\varepsilon '$. 
\einspr

Now we look at the antipode and we proceed in a similar fashion.
\snl
Again we first need a lemma.

\inspr{3.6}  Lemma \rm 
If $\sum (\iota \ot \vp )(\co  (a_i)(1\ot b_i))=0$, then
also $\sum (\iota \ot \vp )((1\ot a_i)\co (b_i))=0$.

\snl\bf Proof: \rm 
As in the proof of Lemma 3.3 we write
$$x\ot c=T_1((\iota \ot \iota \ot \vp )(\co_{13}(a)\co_{23}(b)(1\ot c\ot 1)))
$$
where $x=(\iota \ot \vp )(\co (a)(1\ot b))$. 
Now we apply $\varepsilon \ot \iota $ to 
$(\iota \ot \iota \ot \vp )(\co_{13}(a)
\co_{23}(b)(1\ot c\ot 1))$ and we get $(\iota \ot \vp )((1\ot a)\co (b)(c\ot
 1))$.
\snl
So, if we apply all this with $a_i, b_i$, take sums and use the
injectivity of $T_1$, we get the result. 
\einspr

Now, we can define the antipode.

\inspr{3.7}  Definition \rm 
Define a linear map $S: A\to A$ by
$$
S((\iota \ot \vp )(\co (a)(1\ot b)))=(\iota \ot \vp )((1\ot a)\co (b)).
$$
\einspr
Again, as we see in the proof of the lemma, we have that
$S(x)c=\s \varepsilon (p_i)q_i$, where $x\ot c=\sum \co (p_i)(1\ot q_i)$.
That is the formula for $S$ used in Proposition 3.1. It 
shows that the antipode does not depend on the choice of $\vp $.
\snl
Now we prove the main properties of the antipode. Again there are
different ways to do this.

\inspr{3.8} Proposition \rm 
We have $S(ab)=S(b)S(a)$ for all $a, b\in A$.

\snl\bf Proof: \rm 
Take $a, b, c\in A$ and use the surjectivity of the
map $T_1$. As in the proof of Proposition 3.5 we get
$$\align
ab\ot c = (a\ot c)(b\ot 1)
&= \sum \co (p_i)(b\ot q_i)\\
&=\sum \co (p_ir_{ij})(1\ot s_{ij}).
\endalign$$
Then
$$\align
S(ab)c &= \sum \varepsilon (p_i)\varepsilon (r_{ij})s_{ij}\\
&=\sum \varepsilon (p_i)S(b)q_{i}=S(b)S(a)c.
\endalign$$
We cancel $c$ and get the result. 

\inspr{3.9} Proposition \rm 
$(S\ot S)\co (a)=\si \co (S(a))$ for all $a\in A$.

\snl\bf Proof: \rm 
We have to give a correct meaning to the above
formula. This can be done by multiplying with elements in $A$
or by extending the maps $S\ot S$ and $\si $ to $M(A\ot A)$.
We will not be very strict however in the following argument
because it is quite obvious how to make it precise.
\snl
Take $a, b\in A$ and define $x=(\iota \ot \vp )(\co (a)(1\ot b))$. Then
 $S(x)=(\iota \ot \vp )((1\ot a)\co (b))$ and
$$\align
\co (S(x))&=(\iota \ot \iota \ot \vp )((1\ot 1\ot a)\co^{(2)}(b))\\
&=(\iota \ot S)(\iota \ot \iota \ot \vp )(\co_{23}(a)\co_{13}(b))\\
&=\si (S\ot \iota )(\iota \ot \iota \ot \vp )(\co_{13}(a)\co_{23}(b))\\
&=\si (S\ot S)(\iota \ot \iota \ot \vp )(\co^{(2)}(a)(1\ot 1\ot b))\\
&=\si (S\ot S)\co (x).
\endalign$$
\einspr

Finally we have the following.

\iinspr{3.10}  Proposition \rm  
For all $a\in A$ we have
$$\align
m(\iota \ot S)\co (a) &=\varepsilon (a)1\\
m(S\ot \iota )\co (a) &=\varepsilon(a)1.
\endalign$$
Moreover, if $R$ is any linear map from $A$ to $A$, satisfying one
of the above formulas, then $R=S$.

\snl\bf Proof: \rm 
First remark that one has to give a precise meaning to the above formulas. For the first one this is done by multiplying with an element in $A$ from the left. For the second one, it is necessary to multiply with an element from the right. Also in the proof below, one makes things precise in this way.
\snl
Take $a, b\in A$. Then
$$\align
(\iota \ot S)\co ((\iota \ot \vp )(\co (a)(1\ot b)))
&=(\iota \ot S)(\iota \ot \iota \ot \vp )(\co^{(2)}(a)(1\ot 1\ot b))\\
&=(\iota \ot \iota \ot \vp )(\co_{13}(a)\co_{23}(b)).
\endalign$$
So, if we apply multiplication we get
$$
m(\iota \ot S)\co ((\iota \ot \vp )(\co (a)(1\ot b)))
=(\iota \ot \vp )\co (a)\co (b)=\vp (ab)1$$
and because $\varepsilon (\iota \ot \vp )(\co (a)(1\ot b))=\vp (ab)$, we get
the first formula.
\snl
By the symmetry of our system, we also have a linear map
$S': A\to A$ satisfying
$$m(S'\ot \iota )\co (a)=\varepsilon (a)1$$
for all $a\in A$.
Now, let $R$ be any linear map from $A$ to $A$ satisfying
$m(R\ot \iota )\co (a)=\varepsilon(a)1$ for all $a$. Then, using the Sweedler notation,
$$
\sum R(a_{(1)})a_{(2)}S(a_{(3)})=\sum  \varepsilon (a_{(1)})S(a_{(2)})=S(a)
$$
and also
$$
\s R(a_{(1)})a_{(2)}S(a_{(3)})=\sum R(a_{(1)})\varepsilon (a_{(2)})=R(a).
$$
We get $R=S$ and in particular $S=S'$. So we have also proven the
second formula as well the uniqueness property. 
\einspr
The last part of this proof is standard but again it has to be done correctly by multiplying with elements in $A$ (because we are working with multiplier Hopf algebras).
\snl
Let us finish this section first with a remark on the $^*$-algebra case. When we have a $^*$-algebra with a comultiplication $\co$, we assume that $\co$ is a $^*$-homomorphism. Usually, we also require the integrals to be positive (see section 2). For the moment, let us not make this requirement. On the other hand, it is always possible to take the integrals self-adjoint (in the sense that $\vp(a^*)=\overline{\vp(a)}$ for all $a$). Then, it is easily seen from the definitions that $\varepsilon$ is a $^*$-homomorphism whereas for the antipode, we get $S(a^*)=(S^{-1}(a))^*$. In fact, these results also follow easily from the uniqueness and the other properties of the counit and the antipode we have proven here. 
\snl
We do not prove the other known properties of multiplier Hopf algebras with integrals (like the uniqueness of the integrals, the weak K.M.S. property, ...). We refer to the original paper [VD5] for these results.

\nl\nl

\bf 4. The case of a cointegral \rm
\nl
In their original paper, Larson and Sweedler considered a pair of
an algebra $A$ and a comultiplication $\co$ with the existence of
a left cointegral. For a comultiplication with a counit $\varepsilon$, a
left cointegral is defined as a non-zero element $h\in A$
satisfying $ah=\varepsilon(a)h$ for all $a\in A$.
\snl
We would like to start without the assumption of the existence of
a counit, as we did also in Section 2 where we proved the main
theorem. Therefore, we will start with another definition.
\snl
Again, $A$ is an algebra over $\Bbb C$ with a non-degenerate
product and $\co : A\to M(A\ot A)$ is a comultiplication on $A$ (Definition 1.1). We will assume right away that $\co $ is regular (Definition 1.2).

\inspr{4.1}  Definition \rm 
A non-zero element $h\in A$ is called a
{\it left cointegral} if $\co(a)(1\ot h)=a\ot h$ for all $a\in A$.
\einspr

If there is a counit $\varepsilon$, then this condition is equivalent
with  $ah=\varepsilon(a)h$ for all $a\in A$. In fact, it is not so hard
to show the existence of a counit in the following sense.

\inspr{4.2} Proposition \rm 
Assume that $h$ is a left cointegral as in Definition 4.1.
If the right leg of $\co $ is all of $A$,
then there is a homomorphism $\varepsilon : A\to \Bbb C$ such that $ah=\varepsilon(a)h$ and $(\iota \ot \varepsilon)\co(a)=a$ for all $a\in A$.

\snl\bf Proof: \rm 
From $\co (a)(1\ot h)=a\ot h$ it follows that $bh\in
\Bbb C h$ for all elements $b$ in the right leg of $\co(a)$.
So, by our assumption, $Ah=\Bbb C h$. Hence, there is a linear
map $\varepsilon : A\to \Bbb C$ defined by $ah=\varepsilon(a)h$.
\snl
We have $\varepsilon(ab)h=abh=\varepsilon(b)ah=\varepsilon(a)\varepsilon(b)h$ and so $\varepsilon(ab)=\varepsilon(a)\varepsilon(b)$ for all $a, b\in A$. Moreover
$$ a\ot h=\co (a)(1\ot h)=(\iota \ot \varepsilon )\co (a)\ot h$$
so that $(\iota \ot \varepsilon)\co(a)=a$ for all $a$. 
\einspr

We also have the following.

\inspr{4.3} Proposition \rm 
Again assume that the right leg of $\co $
is all of $A$ and let $\varepsilon: A\to \Bbb C$ be the homomorphism
constructed in the previous proposition. We also have
$(\varepsilon \ot \iota )\co(a)=a$ for all $a\in A$ and $\co(a)(h\ot 1)=h\ot a$.

\snl\bf Proof: \rm 
For all $a\in A$ we have
$$
(\iota \ot \varepsilon \ot \iota )\co^{(2)}(a)=\co(a).
$$
As the right leg of $\co $ is assumed to be all of $A$, this
formula implies that $(\varepsilon \ot \iota)\co (b)=b$ for all $b\in A$.
\snl
Then $\co(a)(h\ot 1)=h\ot (\varepsilon \ot 1)\co(a)=h\ot a$ for all
$a\in A$. 
\einspr

Observe that the formula $\co(a)(1\ot h)=a\ot h$ implies that $a$
belongs to the left leg of $\co (a)$. So the left leg of $\co $ is
all of $A$. If we now also assume that the right leg of $\co $ is
all of $A$, then $\co $ is full (Definition 1.8). So, for a full comultiplication
the equality $\co(a)(1\ot h)=a\ot h$ for all $a$ implies that
also $\co(a)(h\ot 1)=h\ot a$. We also have a counit $\varepsilon$ and $h$
is a cointegral in the original sense.
\snl
Let us mention here that the situation above is very similar to the one encountered in some notes that were written in 1993 (see [VD1]). We have been inspired by some of the techniques used in these notes.
\nl
In what follows, we will now assume that the comultiplication is full.
\snl
As in the case of integrals (Section 2), we can use the
cointegrals to show that the four maps $T_1, T_2, T'_1, T'_2$ are
bijective. We need the cointegrals to be 'faithful':

\inspr{4.4} Definition \rm 
A cointegral $h$ is called {\it faithful} if both
the left and the right leg of $\co(h)$ are all of $A$.
\einspr

Observe that this notion is in accordance with the notion of
faithfulness for a linear functional (Definition 2.2) when we
consider $h$ as a linear functional on the (reduced) dual. 
In the finite-dimensional case, one condition is sufficient. If e.g.\ the left leg of $\co(h)$ is all of $A$, then the same is true for the right leg. Compare with the remark following Definition 2.2.
\snl
Also remark that, when there exists a faitful left cointegral, it is automatic that the comultiplication is full.
\snl
Now, we have the following analogues of Lemma 2.4 and Lemma 2.7.

\inspr{4.5} Lemma \rm 
If $h$ is a left cointegral such that the left
leg of $\co(h)$ is all of $A$, then the map $T_2$ is injective.

\snl\bf Proof: \rm 
Take $a\in A$, apply $\co \ot \iota $ to the equation
$h\ot a=\co(a)(h\ot 1)$ and multiply with $b\in A$ to get
$$\align
(b\ot 1\ot a)(\co(h)\ot 1)
&=(b\ot 1\ot 1)\co^{(2)}(a)(\co(h)\ot 1)\\
&=(\iota \ot \co)((b\ot 1)\co(a))(\co(h)\ot 1).
\endalign$$
So, if $\sum (b_i\ot \iota )\co (a_i)=0$, we get
$$\sum (b_i\ot 1\ot a_i)(\co (h)\ot 1)=0.$$
By assumption, the left leg of $\co (h)$ is all of $A$, so
$\sum b_ic\ot a_i=0$
for all $c\in A$. Therefore $\sum b_i\ot a_i=0$. This proves the
result. 
\einspr

Doing the same thing for the other leg, and on the other side
(flip $\co$, flip the product, or flip both), we get:

\inspr{4.6} Proposition \rm 
If there exist a faithful left integral
and a faithful right integral, then the four maps
$$\matrix
&a\ot b\mapsto \co (a)(1\ot b) &\quad \quad 
&a\ot b\mapsto \co (a)(b\ot 1)\\
&a\ot b\mapsto (a\ot 1)\co(b) &\quad \quad 
&a\ot b\mapsto (1\ot a)\co(b)
\endmatrix$$
are injective.
\einspr

What about surjectivity?

\inspr{4.7} Lemma \rm  
If the left leg of $\co(h)$ is all of $A$, then $T_1$ is surjective.

\snl\bf Proof: \rm
Now apply $\iota \ot \co $ to $a\ot h=\co (a)(1\ot h)$ and
multiply with $1\ot 1\ot b$ to get
$$a\ot (\co (h)(1\ot b))=\co^{(2)}(a)(1\ot (\co (h)(1\ot b))),$$
write $\co(h)(1\ot b)=\sum u_i\ot v_i$ and apply
a linear functional $\om $ (on the third leg of this equation) to get
$$a\ot x=\sum \co (p_i)(1\ot u_i)$$
where $x=(\iota \ot \om )(\co (h)(1\ot b))$ and
$p_i=(\iota \ot \om )(\co (a)(1\ot v_i))$.
This proves the lemma. 
\einspr

Similarly, we find:

\inspr{4.8} Proposition \rm 
If there exist a faithful left integral
and a faithful right integral then the four maps
$$\matrix
&a\ot b\mapsto \co(a)(1\ot b)  &\quad \quad 
&a\ot b\mapsto \co(a)(b\ot 1)\\
&a\ot b\mapsto (a\ot 1)\co(b)  &\quad \quad 
&a\ot b\mapsto (1\ot a)\co(b)
\endmatrix$$
are all surjective.
\einspr
So, we get the following generalization of the Larson-Sweedler result.

\inspr{4.9} Theorem \rm 
Let $A$ be an algebra over $\Bbb C$ with
a non-degenerate product. Let $\co$ be a regular
comultiplication  on $A$. Assume that there are faithful left and right cointegrals. Then $A$ is a regular multiplier Hopf algebra with integrals (of discrete type).

\snl\bf Proof: \rm
By Proposition 4.6 and 4.8, all the four maps
$T_1, T_2, T'_1, T'_2$ are bijective. Hence we have a regular
multiplier Hopf lgebra. Because there exist cointegrals, there
also exist integrals (see [VD3] and [VD-Z1). 
\einspr

When we look at the special case of a finite-dimensional
algebra, we recover (a slightly stronger form of) the original version of the Larson-Sweedler theorem:

\iinspr{4.10} Theorem \rm 
Let $A$ be a finite-dimensional algebra with $1$. Let
$\co $ be a comultiplication on $A$.
Assume that  there is a faithful left cointegral. Then $A$ is a Hopf algebra.

\snl\bf Proof: \rm 
By Lemma 4.5 the map $T_2$ is injective and
so also bijective as $A$ is finite-dimensional. By Lemma 4.7 the
map $T_1$ is surjective and so also bijective. Therefore
$A$ is a Hopf algebra. 
\einspr

Recall that the existence of a faithful left cointegral (in the sense of Definition 4.1) gives us the existence of a counit. The cointegral is also a cointegral in the usual sense. In the original form of the theorem of Larson and Sweedler, the existence of a counit is part of the assumptions. It is also interesting to compare Theorem 4.10 with Theorem 2.9.
\snl
To finish this section, let us indicate how it is possible to develop the theory, starting from the assumptions in Theorem 4.9. Because this case is less general than the case studied in Section 2, we will not do this in detail as we have done in Section 3. We will just briefly give some indications.
\snl
One method would be to start with the assumptions in Theorem 4.9, use the result and pass immediately to the dual. Then one can proceed as in Section 3. A more direct approach would be as follows.
\snl
Let $A$ be an algebra over $\Bbb C$ with a non-degenerate product and let 
$\co$ be a regular coproduct on $A$. Assume as in Theorem 4.9 that there is a faithful left cointegral and a faithful right cointegral. We have seen in Proposition 4.2 and Proposition 4.3 how to find the counit and its properties. What about the antipode? The following property (which is well-known) is dual to the formula in Proposition 3.2 in the case of integrals and again, it is crucial for this point of view. Remark that it was also used in [VD1] to construct the antipode.

\iinspr{4.11} Proposition \rm
Let $(A,\co)$ be a multiplier Hopf algebra with a left cointegral $h$. Then 
$$(1\ot a)\co(h)=(S(a)\ot 1)\co(h)$$
for all $a\in A$.

\snl\bf Proof: \rm
Using the Sweedler notation, we get for all $a\in A$ that
$$\align (1\ot a)\co(h)&=(S(a_{(1)})a_{(2)}\ot a_{(3)})\co(h)\\
 &=(S(a_{(1)})\ot 1)\co(a_{(2)}h)\\
 &=(\varepsilon(a_{(2)})S(a_{(1)})\ot 1)\co(h)\\
 &=(S(a)\ot 1)\co(h)
\endalign$$
and this proves the result.
\einspr

Because of the faithfulness of $h$, the left leg of $\co(h)$ is all of $A$. Now, it is not so difficult to show that the result in Proposition 4.11 can be used to define the antipode and prove its properties. We leave it to the reader as an excercise. 
\snl
And just a small remark to finish. The result in Proposition 4.11 can be used to give a simple argument to show that the antipode is automatically bijective in the case of a finite-dimensional Hopf algebra. Indeed, if $S(a)=0$ it follows from the above formula and the fact that the right leg of $\co(h)$ is all of $A$, that $a=0$.
\nl\nl

\bf References \rm
\bigskip

\parindent 0.5 cm

\item{[A]} E.\ Abe: \it Hopf algebras. \rm Cambridge University Press (1977).
\smallskip

\item{[D-VD]} B.\ Drabant \& A.\ Van Daele:
{\it Pairing and the quantum double of multiplier Hopf algebras}.
Algebras and representation theory 4 (2001), 109-132.
\smallskip

\item{[K-P]} G.I.\ Kac \& V.G.\ Palyutkin: {\it Finite ring groups}. Transactions of the Moscow Mathematical Society 15 (1966), 224-261 (in Russian).
\smallskip

\item{[K-V1]} J.\ Kustermans \& S.\ Vaes: \it A simple
definition  for locally compact quantum groups. \rm  C.R. Acad.
Sci., Paris, S\'er. I 328 (10) (1999), 871-876.
\smallskip

\item{[K-V2]} J.\ Kustermans \& S.\ Vaes: \it Locally compact
quantum groups. \rm Ann.\ Sci.\
 Ec.\ Norm.\ Sup.\  33
(2000), 837-934.
\smallskip

\item{[L-S]} R. G.\ Larson, M. E.\ Sweedler: {\it  An associative
 orthogonal bilinear form for Hopf algebras.} Amer. J. Math. 91(1969) 75-93.
\smallskip

\item{[S]} M.E.\ Sweedler: \it Hopf algebras. \rm Mathematical Lecture Note
Series. Benjamin (1969).
\smallskip

\item{[V-VD]} S.\ Vaes \& A.\ Van Daele: {\it Hopf C$^*$-algebras}.
Proc.\ London Math.\ Soc.\ {\bf 82} (2001) 337-384.
\smallskip

\item{[VD1]} A.\ Van Daele: {\it Quasi-Discrete Locally
Compact Quantum Groups}.  Preprint K.U.Leuven (1993), math.OA/0408208.
\smallskip

\item{[VD2]} A.\ Van Daele: {\it Multiplier Hopf algebras.} Trans.\
Amer.\ Math.\ Soc.\ 342 (1994), 917-932.
\smallskip

\item{[VD3]} A.\ Van Daele: {\it Discrete quantum groups.} J.\ of Alg.\
180 (1996), 431-444.
\smallskip

\item{[VD4]} A.\ Van Daele: {\it The Haar measure on finite quantum
groups}. Proc.\ Amer.\ Math.\ Soc.\ 125 (1997), 3489-3500.
\smallskip

\item{[VD5]} A.\ Van Daele: {\it An algebraic framework for group duality.}
Adv. in Math.  140 (1998), 323--366.
\smallskip

\item{[VD-Z1]} A.\ Van Daele \& Y.\ Zhang : {\it Multiplier Hopf
algebras of discrete type}. J.\ of Alg.\ 214 (1999), 400-417.
\smallskip

\item{[VD-Z2]} A.\ Van Daele \& Y.\ Zhang:
{\it A survey on multiplier Hopf algebras.} Proceedings of the
conference in Brussels on Hopf algebras. Hopf Algebras and Quantum
Groups, eds. Caenepeel/Van Oystaeyen (2000), 269-309. Marcel
Dekker (New York).
\smallskip

\end